\documentclass[10pt]{article}
\usepackage[pdftex,marginratio=1:1,width=180mm,height=235mm,top=15mm,headsep=15mm,head=10mm,footskip=10mm]{geometry}
\usepackage{amsmath,amsfonts,amssymb,amsthm,cite,graphicx,caption,subcaption,bm,array,booktabs,multirow,algorithm,authblk,lineno,xcolor,lipsum,epstopdf,color,hyperref}
\usepackage[english]{babel}
\usepackage{mathtools}
\usepackage{todonotes}
\usepackage{url,floatflt}
\title{Spectral Monic Chebyshev Approximation for Higher Order Differential Equations}

\author[1,2]{M. Abdelhakem\footnote{Correspondence, Email: mabdelhakem@yahoo.com}}
\author[1]{Aya Ahmed}
\author[1,2]{M. El-kady}

\affil[1]{\small{Department of Mathematics, Faculty of Science, Helwan University, Cairo, Egypt}}
\affil[2]{\small{Canadian International College, School of Engineering, Cairo, Egypt}}

\usepackage{academicons}
\usepackage{scalerel}
\usepackage{tikz}
\usetikzlibrary{svg.path}
\definecolor{orcidlogocol}{HTML}{A6CE39}
\tikzset{
  orcidlogo/.pic={
    \fill[orcidlogocol] svg{M256,128c0,70.7-57.3,128-128,128C57.3,256,0,198.7,0,128C0,57.3,57.3,0,128,0C198.7,0,256,57.3,256,128z};
    \fill[white] svg{M86.3,186.2H70.9V79.1h15.4v48.4V186.2z}
                 svg{M108.9,79.1h41.6c39.6,0,57,28.3,57,53.6c0,27.5-21.5,53.6-56.8,53.6h-41.8V79.1z M124.3,172.4h24.5c34.9,0,42.9-26.5,42.9-39.7c0-21.5-13.7-39.7-43.7-39.7h-23.7V172.4z}
                 svg{M88.7,56.8c0,5.5-4.5,10.1-10.1,10.1c-5.6,0-10.1-4.6-10.1-10.1c0-5.6,4.5-10.1,10.1-10.1C84.2,46.7,88.7,51.3,88.7,56.8z};
  }
}
\newcommand\orcidicon[1]{\href{https://orcid.org/#1}{\mbox{\scalerel*{
\begin{tikzpicture}[yscale=-1,transform shape]
\pic{orcidlogo};
\end{tikzpicture}
}{|}}}}

\newcommand{\bq }{\begin{equation}}
\newcommand{\eq }{\end{equation}}

\begin{document}
\maketitle

\begin{abstract} 
This paper is focused on performing a new method for solving linear and nonlinear higher-order boundary value problems (HBVPs). This direct numerical method based on spectral method. The trial function of this method is the Monic Chebyshev polynomials (MCPs). This method was relying on derivative of MCPs which explicit in the series expansion. The advantage of this method is solved HBVPs without transforming it to a system of lower-order ordinary differential equations (ODEs). This method supported by examples of HBVPs in wide application. The mentioned examples showed that the proposed method is efficient and accurate.

 \end{abstract}

\\
%\noindent\textbf{AMS Mathematics Subject
%Classification 2010}: 
\[\]

\section{Introduction}
 In the last few decades, the popularity of the spectral methods in approximate computational has been gained. It has proven highly effective for solving many problems in different branches like physics and chemistry \cite{A1,A2}. Spectral methods solved ODEs by expressing these equations in terms of a series of known and smooth functions. The basic concept of spectral methods is to use a set of trial functions, also called basis or expansion approximating functions. Very smooth, globally and orthogonal are considered to be the essential properties of this polynomials.

Many applications solved by spectral methods which provided better results \cite{A3,A4,A5} time-space with subdiffusion and superdiffusion, Abel's integral equations and the multi-dimensional fractional rayleigh-stokes problem in fluid, respectively, have been solved by spectral methods. Worth mentioning, spectral method has been displayed in \cite{A6}, fifth  kind Chebyshev \cite{A7} and is still being used and developed in \cite{A8}.

The superior approximation of the spectral methods depends on the type of the basis function. It involved three types namely the galerkin, tau and collocation (pseudospectral) methods. Authors in \cite{A9} mentioned galerkin method by monic gegenbauer polynomial. K. Issa and R. B. Adeniyi \cite{A10} presented solution of non-Linear ordinary differential equations depends on tau method. A. Napoli and W. M. Abd-Elhameed in \cite{A11} solved initial value problems based on collocation method.

The HBVPs were treated by an enormous number of algorithms using numerical analysis. Many applications are modeled by HBVPs. It appears in a variety of different fields of engineering, technology, applied science, physics and fluid dynamics are described by boundary value problems (BVPs). Authors in \cite{A12} presented application of HBVPs.

Many authors have treated with high odd-order and high even-order BVPs, particularly \cite{A13}. The solutions of second-order BVPs mentioned in \cite{A14} by predictor-corrector method. The fourth-order BVPs solved numerically in many papers by a different methods \cite{A15}, \cite{A16} and \cite{A17}. The sixth-order is mentioned in \cite{A18} by Legendre wavelet collocation method, also solved in \cite{A19}. On the other hand, many authors solved high odd-order BVPs individuals \cite{A20}. The fifth-order BVPs were solved by \cite{A21}.

It should be noted that there are many numbers of the literature for solving HBVPs. In \cite{A22}, F. A. Costabile and A. Napoli solved HBVPs by birkhoff Lagrange collocation methods. Also, the ultraspherical method has been mentioned in \cite{A23}.

The trial function used in this paper is The MCPs. It is defined by dividing Chebyshev polynomial by its leading coefficient. In \cite{A24}, authors constructed MCPs differentiation matrix. It used in many applications as optimal control \cite{A25}.

The outline of this paper is: in section 2, we expose an overview of MCPs and the applicable characteristics required in what follows. The description of the MCPs method (the m-th derivatives of MCPs) is proven, in section 3. Section 4, the algorithm of the method is given. Section 5 exhibited some numerical results which revealed the accuracy and efficiency of our postulated numerical algorithms. In The 6th and last section, few concluding remarks were presented.

\section{Preliminaries}
In this section, we shall present the base formula for MCPs and the relation between Chebyshev polynomials and MCPs.
The Chebyshev polynomials of first kind \cite{A26} are denoted by $"T_{n}(z);\quad n=1,2,3,\dots \quad on\quad z\in[-1,1]"$ , where
\begin{equation}\label{1}
T_{n}(z)=\cos(n \arccos(z))
\end{equation}
The fundamental recurrence relation is:
\begin{equation}\label{2}
T_{n}(z)=2zT_{n-1}(z)-T_{n-2}(z);\quad n\geq 2
\end{equation}
where $T_{0}(z)=1,\quad T_{1}(z)=z.$\\
The Chebyshev polynomials $"T_{n}(z)"$ can be expanded in power series:
\begin{equation}\label{3}
T_{n}(z)=\frac{n}{2} \sum_{k=0}^{\lfloor n/2 \rfloor} \frac{(-1)^{k}(n-k-1)!(2z)^{n-2k}}{(n-2k)!(k)!}
\end{equation}
where $\lfloor n/2 \rfloor$ is the integer part of $n/2$.
This explicit representation allows us to derive many useful formulas respective to the Chebyshev polynomials \cite{A27}
\begin{equation}\label{4}
2T_{n}(z)=\frac{1}{n+1} T^{'}_{n+1}(z)-\frac{1}{n-1} T^{'}_{n-1}(z) 
\end{equation}
The Chebyshev polynomials form a complete orthogonal set on the interval $[-1,1]$ w.r.t. the weighting function $\varpi(z)=(1-z^2 )^{-1/2}$
\begin{equation}
\label{5}
\langle T_{i},T_{j}\rangle=\int_{-1}^{1}T_{i}(z)T_{j}(z)\varpi(z)dz=\left\{
\begin{array}{cc}
0, & i\neq j\\[.2cm]
\frac{\pi}{2}, & i=j\neq 0 \\[.2cm]
\pi, & i=j=0
\end{array}
\right.    
\end{equation}
Consider MCPs of the first kind \cite{A28} denoted by $"Q_{n}(z) \quad ;\quad n=0,1,2,.... \quad on \quad z\in [-1,1]"$ is given by:
\begin{equation}\label{6}
Q_{n}(z)=\left\{
\begin{array}{cc}
1,& n=0\\[.2cm]
2^{1-n}T_{n}(z), & n\geq 1
\end{array}
\right.    
\end{equation}
Using the previous relation, we have:
\begin{equation}
Q_{1}(z)=z,\quad Q_{2}(z)=z^{2}-\frac{1}{2}
\end{equation}
The recursive formulae of MCPs as follows 
\begin{equation}\label{7}
Q_{n}(z)=zQ_{n-1}(z)-\frac{1}{4}Q_{n-2}(z);\quad n\geq 3
\end{equation}
The MC polynomials constitutes an orthogonal basis w.r.t.$\varpi (z)$:
\begin{equation}\label{8}
\langle Q_{i},Q_{j}\rangle=\int_{-1}^{1}Q_{i}(z)Q_{j}(z)\varpi(z)dz=\left\{
\begin{array}{cc}
0, & i\neq j\\[.2cm]
2^{1-2i}\pi, & i=j\neq 0 \\[.2cm]
\pi, & i=j=0
\end{array}
\right.    
\end{equation}

\section{Monic Chebychev polynomials}
In this section, some important relations and properties of the MCPs have presented and introduced. These relations are essential to setup the Monic Chebychev differentiation expansions. 
The next lemma introduced the recurrence relation for MCPs. This recurrence relation is different from the recurrence relation Chebychev polynomials. This illustrates MCPs is not a copy of Chebychev polynomials. Also, this deviation leads us the totally different differentiation expansion rather those in Chebychev polynomials.

\textbf{Lemma 3.1}

Let $Q_{n}$, be a MC of degree n, then
\begin{equation*}
Q_{n}(z)=\frac{1}{n+1} Q^{'}_{n+1}(z)-\frac{1}{4(n-1)} Q^{'}_{n-1}(z);\quad n\geq 2 
\end{equation*}

\textbf{Proof.}

The proof is straight forward from Eq. (\ref{4}) and Eq. (\ref{6}).

{Theorem 3.1}

The m-th derivative of MCPs is expanded as follows:
\begin{equation}
Q^{(m)}_{n}(z)=\sum_{k=0}^{\lfloor (n-m)/2 \rfloor} b_{nk}c^{(m)}z^{n-2k-m}
\end{equation}
where
\begin{equation}\label{12}
b_{nk}=\frac{(-1)^{k}2^{-2k}n\varGamma(n-k)}{\varGamma(k+1) \varGamma(n-2k+1)}
\end{equation}
\begin{equation}
c^{(m)}=\sum_{l=0}^{m-1}(n-2k-l), \quad m\geq 1
\end{equation} 

\textbf{Proof.}

The following formula is important for proving theorem.
From eq. (\ref{3}), by replacing $z$ with $2^{1-n}z$ then MCPs $Q_{n}(z)$ can be expanded in power series:
\begin{equation}\label{14}
Q_{n}(z)=\sum_{k=0}^{\lfloor n/2 \rfloor} b_{nk}z^{n-2k}, \quad n>0
\end{equation}
where $\lfloor n/2 \rfloor$ is the integer part of $n/2$ and $b_{nk}$ as equation (\ref{12}) 
The first derivative of the formula (\ref{14}) is the following:
\begin{equation}
Q^{(1)}_{n}(z)=\sum_{k=0}^{\lfloor n/2 \rfloor} (n-2k)b_{nk}z^{n-2k-1}, \quad n>0
\end{equation}	
Then, using the mathematical induction to prove the theorem:
For $m=1$, we have
\begin{equation}
Q^{(1)}_{n}(z)=\sum_{k=0}^{\lfloor (n-1)/2 \rfloor} (n-2k)b_{nk}z^{n-2k-1}, \quad n>0
\end{equation}
which leads to the same result as the direct first differentiation.
Now, assume that theorem is satisfied at $m=j$
\begin{equation}
Q^{(j)}_{n}(z)=\sum_{k=0}^{\lfloor (n-j)/2 \rfloor}\frac{(-1)^{k}2^{-2k}n\varGamma(n-k)}{\varGamma(k+1) \varGamma(n-2k+1)} z^{n-2k-j}, 
\end{equation}
For $m=j+1$, we have

\begin{eqnarray}\label{8}
Q^{(j+1)}_{n}(z)&=&\frac{d}{dz}Q^{(j)}_{n}(z)=\frac{d}{dz} 
\sum_{k=0}^{\lfloor(n-j)/2\rfloor}\left( \frac{(-1)^{k}2^{-2k}n\varGamma(n-k)}{\varGamma(k+1) \varGamma(n-2k+1)} z^{n-2k-j}\nonumber \right)\\
&=&\sum_{k=0}^{\lfloor \frac{n-(j+1)}{2} \rfloor} \Big(\frac{(-1)^{k}2^{-2k}n\varGamma(n-k) }{\varGamma(k+1) \varGamma(n-2k+1)}(n-2k-j)\, c^{(j)} z^{n-2k-j-1} \Big) \nonumber\\
&=&\sum_{k=0}^{\lfloor \frac{n-(j+1)}{2}\rfloor}\big(\frac{(-1)^{k}2^{-2k}n\varGamma(n-k)}{\varGamma(k+1) \varGamma(n-2k+1)}(n-2k-j) \,c^{(j+1)} z^{n-2k-j-1} \big)
\end{eqnarray}
Hence the theorem is satisfied for every positive integer m.

\section{Proposed method}

For solving the problems in the next section, we need to calculate the $y(z_{i})$ such that $z_{i}$ is Chebyshev Gauss Lobatto points where $z_{i}=\cos(\pi i/N);\quad 0\leq i\leq N$\\
Consider the linear HBVPs:
\begin{eqnarray}\label{19}
f_{m}(z)y^{(m)}(z)&+&f_{m-1}(z)y^{(m-1)}(z)+f_{m-2}(z)y^{(m-2)}(z)\nonumber \\[.1cm]
&+&\dots +f_{0}(z)y(z)=g(z);a\leq x\leq b
\end{eqnarray}
Subject to the boundary conditions:
\begin{equation}\label{20}
\left\{
\begin{array}{l}
y(a)=\alpha_{0},\quad y(b)=\beta_{0},\\[.2cm]
y^{(1)}(a)=\alpha_{1},\quad y^{(1)}(b)=\beta_{1},\\[.2cm]
\vdots\\[.2cm]
y^{(q)}(a)=\alpha_{q},\quad y^{(q)}(b)=\beta_{q}.
\end{array}
\right.    
\end{equation}
such that $f(z)$ and $g(z)$ is a function of $z$ or constant, $\alpha_{0}$ and  $\beta_{0}$ are constants. 
The number of boundary conditions is equal to the order of the problem. 
Consider the non-linear HBVPs as:
\begin{equation}\label{21}
y^{(m)}_{n}=g(z,y,y^{'});\quad a\leq x\leq b
\end{equation}
Such that $g(z,y,y^{'})$ is a non-linear function. For examples, $e^{y^{(m)}(z)}$ and $(y^{(m)}(z))^{n}, n>1 \quad etc \dots $,
Subject to the boundary conditions as in (\ref{20}).
Furthermore, by approximating the unknown function $y(z_{i})$ as
\begin{equation}\label{22}
y(z_{i})=\sum_{i=0}^{N}A_{i}Q_{i}(z_{j}),\quad z_{j}=\cos(\frac{\pi j}{N});0\leq j\leq N
\end{equation}
Where $A_{i}$  denotes constants. Realizing that $Q_{i}(z_{j})$  is square matrix. i.e. the elements of this matrix are the values of Monic Chebyshev polynomials at $z_{j}$
The derivatives can be expressed in differential equation as the following form:
\begin{equation}\label{23}
\frac{d^{m}y}{dz^{m}}=\sum_{i=0}^{N}A_{i}Q_{i}(z)
\end{equation}
Such that m is the order of the derivative. At m=0, the form (\ref{23}) become equivalent to (\ref{22}). 
Substituting from (\ref{23}) into (\ref{19}) yields:
\begin{eqnarray}
f_{m}\sum_{i=0}^{N}A_{i}Q^{(m)}_{i}(z)&+&f_{m-1}\sum_{i=0}^{N}A_{i}Q^{(m-1)}_{i}(z)\nonumber\\
&+&f_{m-2}\sum_{i=0}^{N}A_{i}Q^{(m-2)}_{i}(z)+\dots\nonumber\\
&+&f_{0}\sum_{i=0}^{N}A_{i}Q_{i}(z)=g(z);\quad a\leq x\leq b
\end{eqnarray}
By shifting the interval $[a, b]$ such that, $z$ will be belongs to the interval from -1 to 1 and shifting the boundary conditions (20) as the following:
\begin{equation}\label{25}
\left\{
\begin{array}{l}
\sum_{i=0}^{N}A_{i}Q_{i}(-1)=\alpha_{0},\quad \sum_{i=0}^{N}A_{i}Q_{i}(1)=\beta_{0},\\[.2cm]
\sum_{i=0}^{N}A_{i}Q^{(1)}_{i}(-1)=\alpha_{1},\quad \sum_{i=0}^{N}A_{i}Q^{(1)}_{i}(1)=\beta_{1}\\[.2cm]
\vdots\\[.2cm]
\sum_{i=0}^{N}A_{i}Q^{(q)}_{i}(-1)=\alpha_{q},\quad\sum_{i=0}^{N}A_{i}Q^{(q)}_{i}(1)=\beta_{q}.
\end{array}
\right.    
\end{equation}
Substituting  (\ref{22}), (\ref{23}) in (\ref{19})
\begin{eqnarray} \label{26}
\sum_{i=0}^{N}A_{i}Q^{(m)}_{i}(z_{j})&+&f_{m-1}\sum_{i=0}^{N}A_{i}Q^{(m-1)}_{i}(z_{j})\nonumber\\
&+&f_{m-2}\sum_{i=0}^{N}A_{i}Q^{(m-2)}_{i}(z_{j})+\dots\nonumber \\[.1cm]
&+&f_{0}\sum_{i=0}^{N}A_{i}Q_{i}(z_{j})=g(z_{j});\quad 0\leq j\leq N
\end{eqnarray}
then equations (\ref{26}) construct a system of algebraic equations. Particularly, linear or non-linear system of equations depend on formulation of HBVPs "(\ref{19}) and (\ref{21})" of the problem. Hint, we shall be replacing some rows of $D_{ij}$  by the conditions (\ref{25}). In case of the non-linearity, the non-linear algebraic system may be treated by numerous methods like newton method, secant method and SOR- Steffensen-Newton method for finding solution of systems of nonlinear equations \cite{A29}. $A_{i}$  can be approximated by solving the previous system. Finally we can get the approximate solution expanded depend on MCPs.

\section{Numerical examples}

In this part, the approximation results of linear and nonlinear HBVPs are solved by applying the previous method. The results emphasized that the method satisfied high accuracy and efficiency.

\textbf{Example 1}

Non-linear fourth order BVP:
\begin{equation}\label{27}
16y^{(4)}(z)-6e^{-4y(z)}=-12(1.5+0.5z)^{-4}; \quad -1\leq z \leq 1
\end{equation}
With boundary conditions:
\begin{equation}\label{28}
\left\{
\begin{array}{cc}
y(-1)=0, & y(1)=ln(2),\\[.2cm]
y^{\prime}(-1)=0.5, & y^{\prime}(1)=0.25.\\[.2cm]
\end{array}
\right.    
\end{equation}
with analytical solution:
\begin{equation*}
y(z)=ln(1.5+0.5z).
\end{equation*}
Applying the method, so equation (\ref{27}) and (\ref{28}) can be written as:
\begin{eqnarray}
16\sum_{i=0}^{N}A_{i}Q^{4}_{i}(z_{j})-6e^{\sum_{i=0}^{N}}A_{i}Q_{i}(z_{j})
&=&-12(1.5+0.5z_{j})^{-4}\nonumber \\[.1cm]
&;& 0\leq j \leq N
\end{eqnarray}
Subject to:
\begin{equation}
\left\{
\begin{array}{cc}
\sum_{i=0}^{N}A_{i}Q_{ij}(-1)=0, & \sum_{i=0}^{N}A_{i}Q_{ij}(1)=ln(2),\\[.2cm]
\sum_{i=0}^{N}A_{i}Q^{\prime}_{ij}(-1)=0.5,&\sum_{i=0}^{N}A_{i}Q^{\prime}_{ij}(1)=0.25.\\[.2cm]
\end{array}
\right.    
\end{equation}

\begin{table}
	\caption{The approximate solutions for \textbf{example 1} at N=10.}
	\begin{center}
			\renewcommand{\arraystretch}{1.5}
		\begin{tabular}{|>{\centering\arraybackslash}p{0.7cm}|>{\centering\arraybackslash}p{2cm}|>{\centering\arraybackslash}p{2cm}|>{\centering\arraybackslash}p{2cm}|}
			\hline
			z & Analytical solution & Approximate solution \cite{A30} & Approximate solution \\
			\hline
			-1 & 0 & 0 & 0.0002098641\\
			-0.8 & 0.0953101798 &0.0950147533  &0.0953558832 \\
			-0.6 & 0.1823215568 &0.1814496227  &	0.1823739721 \\
			-0.4 &	0.2623642645 &0.2609546573&0.2624093359  \\ 
			-0.2 & 0.3364722366 &0.3347370220  &0.3364999034\\
			0 &	0.4054651081 & 	0.4036840381&  	0.4054692075 \\
			0.2 & 0.4700036292  &0.4684459279 &  0.4699817600 \\
			0.4 & 0.5306282511 &	0.5294932609 & 	0.5305817961 \\
			0.6 &	0.5877866649 &	0.5871580370  &0.5877209165  \\
			0.8 &	0.6418538862 & 	0.6416636708& 	0.6417780693  \\
			1 &	0.6931471806 &	0.6931471806  &	0.6930743766  \\
			\hline
	
		\end{tabular}
	\end{center}
\end{table}

\begin{table}
	\caption{The AE for \textbf{example 1} at N=10.\label{30}}
	\begin{center}
		\renewcommand{\arraystretch}{1.5}
		\begin{tabular}{|c|c|c|}
			\hline
			z &  AE \cite{A30} &AE (Present method)\\
			\hline
			-1 & 0&  2.098641e-05 \\
			-0.8 &  2.954265e-04 & 4.570342e-05 \\
			-0.6 &  8.719341e-04 & 5.241534e-05 \\
			-0.4 &	 1.4096072e-03&4.507148e-05 \\ 
			-0.2 & 1.7352146e-03 & 2.766681e-05 \\
			0 &	 1.7810699e-03 & 4.099373e-06 \\
			0.2 & 1.5577013e-03 & 2.186930e-05 \\
			0.4  & 1.1349902e-03 &	4.645491e-05 \\
			0.6  & 6.286279e-04 &6.574845e-05 \\
			0.8& 1.902154e-04 &	7.581684e-05 \\
			1 &	 0  &	7.280400e-05 \\
			\hline
			
		\end{tabular}
	\end{center}
\end{table}

\begin{table}
	\caption{The MAE for \textbf{example 1} .\label{31}}
	\begin{center}
			\renewcommand{\arraystretch}{1.5}
		\begin{tabular}{|c|c|}
			\hline
		N & MAE\\
		\hline
		10 & 7.5991e-05 \\
		12 &  3.5789e-07 \\
		14 & 1.2290e-09 \\
	    16 &  3.5426e-12 \\
	    18 & 1.5099e-14 \\
	    20 & 4.4409e-16 \\
			\hline
		\end{tabular}
	\end{center}
\end{table}

As shown in table (\ref{30}), the MAE for the present method is more efficient compared by the \cite{A30}. In table (\ref{31}), the present method gets e-16 as the MAE at N=20 which proved that method is accuracy.
%%%%

\textbf{Example 2}

The following third order BVP:
\begin{eqnarray}
8y^{(3)}(z)
&-&0.5(z+1)y(z))=(0.125(z+1)^{3}\nonumber \\[.1cm]
&-&0.5(z+1)^{2}-2.5(z+1)-3)e^{0.5(z+1)} \nonumber \\[.1cm]
&;& \quad -1\leq z \leq 1
\end{eqnarray}
Subject to:
\begin{equation}
y(-1)=y(1)=0, \quad y^{(1)}(-1)=0.5. 
\end{equation}
With analytical solution
\begin{equation*}
y(z)=0.25(1-z^{2})e^{0.5(z+1)}.
\end{equation*}
Appling the same algorithm as a previous example, the results are summarized as table (\ref{32}):

\begin{table}
	\caption{The MAE for \textbf{example 2} .\label{32}}
	\begin{center}
		\renewcommand{\arraystretch}{1.5}
		\begin{tabular}{|c|c|c|c|c|}
			\hline
			\multirow{2}{*}{N} & \multicolumn{3}{c|}{ \cite{A30} } & \multirow{2}{*}{present method} \\ 
			\cline{2-4}
			&3SHLM1& 3SHLM2& 3SHLM3& \\
			\hline
			6& 1.079e-06 & 9.13e-08 & 8.33e-08 & 2.563e-04 \\ 
			9& 1.210e-7 & 6.80e-09  & 3.30e-09 & 2.809e-08 \\ 
			10& --- & --- & --- & 9.966e-10 \\ 
			12& 2.770e-8 & 1.10e-09 & 9.53e-10 & 8.312e-13 \\ 
			14& --- & --- & --- & 4.441e-16 \\ 
			15& 8.900e-09 & 2.00e-10 & 1.13e-10 & 1.110e-16 \\
			\hline 
		\end{tabular}
	\end{center}
\end{table}

In table (\ref{32}), we notice that, at N=12, the MAE is 8.312e-13 for our method but the MAE is 9.53e-10 for \cite{A31}. The presented algorithm gets the double precision "e-16" as a MAE at N=15. That compression emphases the efficiency and accuracy of the presented method.

\textbf{Example 3}

The following sixth order BVP:
\begin{equation}
y^{(6)}(z)+y(z)=12z\cos(z)+30\sin(z); \quad -1\leq z \leq 1
\end{equation}
Subject to:
\begin{equation}
\left\{
\begin{array}{c}
y(\pm 1)=0, \quad y^{\prime}(\pm 1)=2\sin(1),\\[.2cm]
y^{''}(-1)=-y^{''}(1)=-4\cos(1)-2\sin(1).\\[.2cm]
\end{array}
\right.    
\end{equation}
With analytical solution:
\begin{equation*}
y(z)=(z^{2}-1)\sin(z).
\end{equation*}
Using the same method as a previous example, the results are summarized as the following:

\begin{table}
	\caption{The MAE for \textbf{example 3} .\label{33}}
	\begin{center}
		\renewcommand{\arraystretch}{1.5}
		\begin{tabular}{|c|c|c|c|}
			\hline
			N & \cite{A32} & \cite{A32} & Present method \\
			\hline
			12 & 1.866e-8 & 1.883e-8 & 3.163e-09 \\
			15 & --- & --- & 4.235e-14 \\
			16 & 2.384e-13 & 2.394e-13 & 1.920e-14 \\
			17 & --- & --- & 4.996e-16 \\
			20 & 2.797e-16 & 2.797e-16 & 3.330e-16 \\
			24 & 2.797e-16 & 3.674e-16 & 6.106e-16 \\
			\hline
		\end{tabular}
	\end{center}
\end{table}
In table (\ref{33}), we notice that, at N=12,  the MAE is 3.1631e-09 for our method but the MAE is 1.866e-8 for \cite{A32}. The presented method gets the double precision "e-16" as a MAE at N=17, while method in \cite{A32} gets it at N=20.  That compression emphases the efficiency and accuracy of the presented method.

\textbf{Example 4}

The ninth order BVP: 
\begin{equation}
512y^{(9)}(z)-y(z)=-9e^{0.5(z+1)}; \quad -1\leq z \leq 1
\end{equation}
Subject to:
\begin{equation}
\left\{
\begin{array}{cc}
y(-1)=1, & y(1)=0, \\[.2cm]
y^{\prime}(-1)=0, & y^{\prime}(1)=-0.5e,\\[.2cm]
y^{''}(-1)=-0.25, & y^{''}(1)=-0.5e,\\[.2cm]
y^{'''}(-1)=-0.25, & y^{'''}(1)=-0.375e\\[.2cm] y^{(4)}(-1)=-0.1875.
\end{array}
\right.    
\end{equation}
with analytical solution:
\begin{equation*}
y(z)=(0.5-0.5z)e^{0.5(z+1)}.
\end{equation*}

\begin{table}
	\caption{The Aproximate solution for \textbf{example 4} at N=10.}
	
	\begin{center}
		\renewcommand{\arraystretch}{1.5}
	\begin{tabular}{|>{\centering\arraybackslash}p{0.7cm}|>{\centering\arraybackslash}p{2cm}|>{\centering\arraybackslash}p{2cm}|>{\centering\arraybackslash}p{2cm}|}
			\hline
			z & Analytical  solution&Approximate  solution \cite{A33}& 	Approximate solution  \\
			\hline
			-1 &1.000000&	1.000000&1.000000	\\
			-0.8&0.994654	& 0.994654&			0.994654	\\
			-0.6&0.977122	&0.977122  &			0.977122 \\
			-0.4&	0.944901&	0.944901	&	0.944901	 \\
			-0.2&0.895095	&	0.895095&	0.895095\\
			0&0.824361&0.824361	&	0.824361	\\
			0.2&	0.728848&	0.728848&0.728848	\\
			0.4&0.604126	&0.604126&0.604126	 \\
			0.6&0.445108	&0.445108	&	0.445108	 \\
			0.8&	0.245960&	0.245960	&0.245960 \\
			1&0.000000	&-6.940873e-09 &	-1.765061e-16 \\
			\hline
	
		\end{tabular}
	\end{center}
\end{table}

\begin{table}
	\caption{The AE for \textbf{example 4} at N=10.\label{34}}
	
	\begin{center}
		\renewcommand{\arraystretch}{1.5}
		\begin{tabular}{|c|c|c|}
			\hline
			z & AE \cite{A33}& AE (present method) \\
			\hline
			-1 &0 &1.110223e-16\\
			-0.8 &2.289423e-10&	6.585842e-12\\
			-0.6 & 4.623567e-09&	1.205764e-10\\
			-0.4 & 2.081388e-08&	4.923047e-10\\
			-0.2 &	4.783428e-08&	1.029267e-09\\
			0 &	7.122960e-08&	1.397438e-09\\
			0.2 &7.339381e-08&	1.321320e-09\\
			0.4 &	4.961130e-08&	8.449000e-10 \\
			0.6 & 1.629848e-08&	3.074025e-10\\
			0.8 &	1.814813e-08 & 	3.313241e-11\\
			1 & 6.940873e-08& 	1.765061e-16\\
			\hline
			
		\end{tabular}
	\end{center}
\end{table}

\begin{table}
	\caption{The MAE for \textbf{example 4} .\label{35}}
	\begin{center}
		\renewcommand{\arraystretch}{1.5}
		\begin{tabular}{|c|c|}
			\hline
			N & MAE\\
			\hline
			10 & 1.3974e-09 \\
			12 & 6.5798e-12 \\
			14 & 1.3212e-14 \\
			16 & 4.4409e-16 \\
			18 & 3.3307e-16 \\
			20 & 4.4409e-16 \\
			\hline
		\end{tabular}
	\end{center}
\end{table}

Table (\ref{34}) shows the MAE for the present method is more efficient compared by the \cite{A33}. Both tables (\ref{34}) and (\ref{35}) proved that the presented method is efficient and accurate.

\textbf{Example 5}

The following linear 12th-order boundary value problem:
\begin{eqnarray}
2096y^{(12)}(z)&+&0.5(z+1)y(z)=-(120+11.5(z+1) \nonumber\\
&+&0.125(z+1)^{3})e^{0.5(z+1)};\quad -1\leq z \leq 1
\end{eqnarray}

Subject to:
\begin{equation}
\left\{
\begin{array}{cc}
y(\pm 1)=0, & y^{\prime}(-1)=0.5,\\[.2cm]
y^{\prime}(1)=-0.5e, & y^{''}(-1)=0, \\[.2cm]
y^{''}(1)=-e, & y^{'''}(-1)=-0.375, \\[.2cm] y^{'''}(1)=-1.125e, & y^{(4)}(-1)=-0.5,\\[.2cm]
y^{(4)}(1)=-e, & y^{(5)}(-1)=-0.46875,\\[.2cm] y^{(5)}(1)=-0.78125e.
\end{array}
\right.    
\end{equation}
With exact solution:
\begin{equation*}
y(z)=0.25(1-z^{2})e^{0.5(z+1)}.
\end{equation*}

\begin{table}
	\caption{The MAE for \textbf{example 5}.\label{36}}
	\begin{center}
		\renewcommand{\arraystretch}{1.5}
		\begin{tabular}{|c|c|c|}
			\hline
			N &	\cite{A34} & Present method\\[0.5ex]
			\hline
			12 & 3.122E-11 & 1.711e-12 \\
			14 & 3.104E-12 & 9.103e-15 \\
			16 & 4.324E-13 & 3.330e-16 \\
			18 & 1.279E-13 & 1.110e-16 \\
			20 & 2.312E-14 & 1.110e-16 \\
			22 & 4.153E-15 & 1.665e-16 \\
			\hline
		\end{tabular}
	\end{center}
\end{table}

Table (\ref{36}) shows the results of MAE between the present method and the method in \cite{A34}. Since the MAE at N=16 is 3.3307e-16 while the method in \cite{A34} gets to 4.153e-15 at N=22 which proven that our present method is more efficient and accurate.

\section*{Conclusion}

We have shown a Galerkin method expanded depend on the MCPs to solve HBVPs. The developed method has been utilized to solve HBVPs without reducing to system of lower order differential equations. The essential advantages of the presented method are the simplicity and easy in apply. Four examples linear and nonlinear of HBVPs have been solved. The obtained results of numerical examples showed that the presented method is an accurate and efficient for the HBVPs.
\section*{Orcid}

M. Abdelhakem \orcidicon{0000-0001-7085-1685}\,\url{https://orcid.org/0000-0001-7085-1685}

\emergencystretch=\hsize

\begin{center}
	\rule{6 cm}{0.02 cm}
\end{center}

\begin{floatingfigure}[h]{3.5cm}
	\centering
	\includegraphics[width=3cm,height=3.4cm]{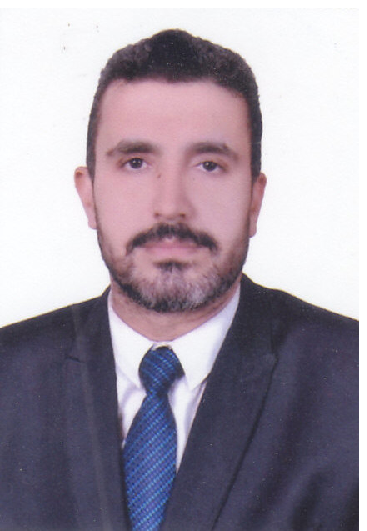}
\end{floatingfigure}
{\bf \hypertarget{author1}{Mohamed Abdelhakem}} "mabdelhakem@yahoo.com" graduated from Helwan University, Cairo, Egypt, in 1999. He received the M. Sc. degree in Computer Science, from Helwan University, Cairo, Egypt, in 2003, and the Ph. D. degree in pure mathematics from Helwan
University, Cairo, Egypt, in 2011. He is the vice Chairman and founder of Helwan School of Numerical Analysis in Egypt "HSNAE".

\end{document}